\documentclass[reqno,12pt]{amsart}

\usepackage{epsf}
\usepackage{graphics}
\usepackage{amssymb}
\usepackage{amsmath}

\date{}

\theoremstyle{plain}
\newtheorem*{main}{Main Theorem}
\newtheorem{theorem}{Theorem}
\newtheorem{corollary}{Corollary}

\newtheorem{proposition}{Proposition}

\theoremstyle{definition}

\theoremstyle{remark}
\newtheorem*{example}{Example}

\newtheorem{question}{Question}
\newtheorem*{remark}{Remark}

\def\N{{\mathbb N}}
\def\Z{{\mathbb Z}}
\def\Q{{\mathbb Q}}
\def\R{{\mathbb R}}

\title{Asymptotic Rasmussen Invariant} 

\author{Sebastian Baader}

\begin{document}

\begin{abstract} We use simple properties of the Rasmussen invariant of knots to study its asymptotic behaviour on the orbits of a smooth volume preserving vector field on a compact domain in the 3-space. A comparison with the asymptotic signature allows us to prove that asymptotic knots are non-alternating, in general. Further we show that the Rasmussen invariant defines a quasi-morphism on the braid groups and derive estimates for the stable commutator and torsion lengths of alternating braids.
\end{abstract}

\maketitle

\section{Introduction}

The Arnold invariant of a smooth volume-preserving vector field on a closed homology 3-sphere measures how two orbits are linked asymptotically (see~\cite{AK}). An asymptotic linking number does also exist for a smooth volume-preserving vector field $X$ on a compact domain $G \subset \R^3$ with smooth boundary, provided $X$ is tangent to the boundary $\partial G$. In~\cite{GG1}, Gambaudo and Ghys proved that the asymptotic linking number can be determined by looking at a single orbit only. More precisely, they proved the existence of an asymptotic signature invariant which coincides with the Arnold invariant divided by two, at least if the flow of $X$ is ergodic.

In this paper we prove the existence of an asymptotic Rasmussen invariant and compare it with the asymptotic signature. As an application, we show that asymptotic knots are non-alternating, as soon as their asymptotic Rasmussen invariant is non-zero. In order to state our main theorem, we have to describe how pieces of orbits can be turned into knots: let $x \in G$ be a non-periodic, non-singular point. For a fixed time $T>0$, we define $K(T,x) \subset \R^3$ to be the piece of orbit from $x$ to $\Phi_X^T(x)$, followed by the geodesic segment $\gamma$ joining $\Phi_X^T(x)$ to $x$ ($\gamma$ need not be contained in $G$). A careful exposition on the closure of pieces of orbits can be found in~\cite{V}. For almost all $x \in G$ and $T>0$, $K(T,x)$ is an embedded curve, i.e. a knot. As usual, we denote the Rasmussen invariant and the signature of a knot $K$ by $s(K)$ and $\sigma(K)$, respectively.

\begin{main} Let $X$ be a smooth vector field on a compact domain $G \subset \R^3$, tangent to the boundary $\partial G$, with hyperbolic singularities only, i.e. linear singularities corresponding to critical points of index 1 or 2 of a Morse function on $\R^3$. If $\mu$ is an $X$-invariant probability measure which does not charge the periodic orbits and singular points of $X$, then the limit
$$s(X,x):=\lim_{T \to \infty} \frac{1}{T^2} s(K(T,x))$$
exists for almost all $x \in G$ (with respect to $\mu$) and coincides with
$$2\sigma(X,x):=2\lim_{T \to \infty} \frac{1}{T^2} \sigma (K(T,x)).$$ 
\end{main}

\begin{remark} As Gambaudo and Ghys observed in~\cite{GG1}, this quantity coincides with the Arnold invariant of $X$, if the flow of $X$ is ergodic with respect to $\mu$, i.e. if every measurable function which is invariant under the flow of $X$ is constant almost everywhere.
\end{remark}

\begin{corollary} If, for a point $x \in G$, the limit $s(X,x)$ exists and is non-zero, then there exists a positive constant $S \in \R$ such that the knot $K(T,x)$ is non-alternating, for all $T \geq S$.
\end{corollary}

As an example, let $X$ be the constant vector field $(1, \omega)$ on $S^1 \times S^1=\R^2/\Z^2$. We can easily extend $X$ to a non-vanishing vector field on the full torus $V=S^1 \times D^2$, which we may view as a submanifold of $\R^3$. If we choose an irrational slope $\omega \in \R-\Q$, then the orbit starting at a point $x \in S^1 \times S^1$ is non-periodic and shows a non-alternating behaviour, as time increases. Its asymptotic Rasmussen invariant is $s(X,x)=2\sigma(X,x)=\frac{1}{4\pi^2} \omega$ (see~\cite{GG1}, p.50).

\begin{remark} The assumption in the corollary that $s(X,x)$ be non-zero is absolutely essential, as shows the example of the constant vector field $(1,0)$ on $S^1 \times D^2=\R/\Z \times D^2$, whose orbits are all periodic and unknotted.
\end{remark} 

The proof of the main theorem is heavily based upon Gambaudo and Ghys' work (\cite{GG1}). We present it in the next section. The last three sections are devoted to the study of the Rasmussen invariant as a quasi-morphism on the braid groups. They can be read independently.

\section{The asymptotic Rasmussen invariant}

The Rasmussen invariant $s$ of knots was constructed from the Khovanov complex of knots in~\cite{Ra}. Among various interesting properties of the Rasmussen invariant, the following three are of special interest to us:
\begin{equation}
s(K) \geq 1+w(D)-o(D),
\label{property1}
\end{equation}
where $w(D)$ and $o(K)$ stand for the writhe (i.e. the algebraic crossing number) and the number of Seifert circles of a diagram $D$ of $K$, respectively,
\begin{equation}
s(K) \leq 2g_*(K),
\label{property2}
\end{equation}
where $g_*(K)$ denotes the 4-genus of a knot $K$,
\begin{equation}
s(K)=\sigma(K),
\label{property3}
\end{equation}
for all alternating knots $K$.

The first inequality was proved by Shumakovitch in~\cite{S}, the other properties appear in Rasmussen's original paper. 
Applying the first inequality to the mirror image $\bar{K}$ of a knot $K$, we obtain
$$s(\bar{K}) \geq 1+w(\bar{D})-o(\bar{D})=1-w(D)-o(D).$$
Combining this with the fact that $s(\bar{K})=-s(K)$, we also get an upper bound for $s(K)$, altogether:
\begin{equation}
1+w(D)-o(D) \leq s(K) \leq -1+w(D)+o(D).
\label{property4}
\end{equation}

\begin{proof}[Proof of the main theorem]
According to Gambaudo and Ghys (\cite{GG1}), the complement of the singularities of $X$ can be covered by an enumerable family of flow boxes whose flow time (i.e. the minimal time it takes to pass through a flow box) is bounded from below by a global constant $\lambda>0$. Further they show that for almost all $x \in G$ and $T>0$, the knots $K(T,x)$ have diagrams $\pi_0(K(T,x))$ whose writhe, called $\theta$ there, grows quadratic in $T$:
\begin{equation}
\lim_{T \to \infty} \frac{1}{T^2} w(\pi_0(K(T,x)))=2\sigma(X,x).
\label{writhe}
\end{equation}
In fact, they subdivide the crossings of $\pi_0(K(T,x))$ into three types (\cite{GG1}, p.64): $D_1$, $D_2$, $D_3$. The crossings of type $D_1$ arise from overcrossing flow boxes, as illustrated in figure~1; their number grows quadratic in~$T$. The number of crossings of types $D_2$ and $D_3$ grows subquadratic in~$T$.

\begin{figure}[ht]
\scalebox{1.2}{\raisebox{-0pt}{$\vcenter{\hbox{\epsffile{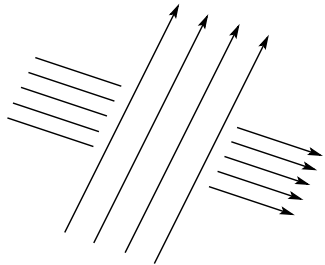}}}$}} 
\caption{two overcrossing flow boxes}
\end{figure}

We shall estimate the number of Seifert circles $o(\pi_0(K(T,x)))$. Every Seifert circle of $\pi_0(K(T,x))$ is adjacent to at least one crossing, and there are at most two Seifert circles meeting at each crossing. Therefore the number of Seifert circles adjacent to a crossing of type $D_2$ or $D_3$ grows subquadratic in $T$. Further, every Seifert circle which is adjacent to a crossing of type $D_1$ must enter at least one flow box (see figure~2). Therefore there are at most $\frac{T}{\lambda}$ such circles, where $\lambda$ is the minimal flow time for flow boxes! Altogether, this shows
$$\lim_{T \to \infty} \frac{1}{T^2} o(\pi_0(K(T,x)))=0.$$
In view of (\ref{property4}) and (\ref{writhe}), this proves the main theorem.
\end{proof}

\begin{figure}[ht]
\scalebox{1.2}{\raisebox{-0pt}{$\vcenter{\hbox{\epsffile{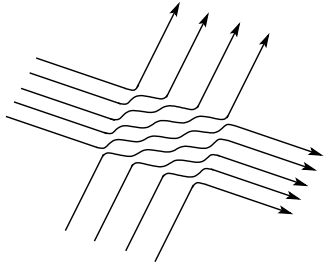}}}$}} 
\caption{pieces of Seifert circles}
\end{figure}

The corollary in turn is an immediate consequence of the main theorem and the fact that $s(K)=\sigma(K)$, for all alternating knots $K$. Indeed, if $s(X,x) \neq 0$, then $s(X,x)=2\sigma(X,x) \neq \sigma(X,x)$, hence there exists a constant $S \geq 0$, such that $s(K(T,x)) \neq \sigma(K(T,x))$, for all $T \geq S$.

\begin{remark} The proof of the main theorem works for all knot invariants $I$ that satisfy the inequality
$I(K) \geq 1+w(D)-o(D)$ and the equation $I(\bar{K})=-I(K)$. This is notably the case for the invariant $2\tau$ coming from the knot Floer homology (\cite{OS}).
\end{remark}

\section{Stable lengths in braid groups}

Given an element $\beta$ of the braid group $B_n$, we may ask what is the minimal number of positive and negative standard generators needed to factorize $\beta$. This number is called the braid length of $\beta$. For instance, the braid length of $\sigma_1 \sigma_2^{-1} \sigma_1 \sigma_2^{-1} \in B_3$ equals $4$. An efficient lower bound for the braid length of positive braids is provided by the linking number $lk: B_n \to \Z$, a homomorphism defined by $lk(\sigma_i)=1$, for all standard generators $\sigma_i \in B_n$. In an analogous way, we may consider the $S$-length $l_S: G \to \N$ for any group $G$ with a generating set $S \subset G$. The \emph{stable $S$-length} $||.||_S: G \to \R$ is then defined by
$$||g||_S=\lim_{n \to \infty} \frac{1}{n} \; l_S(g^n).$$
In particular, the \emph{stable commutator length} in $G$ is the stable length with respect to the set $C$ of commutators in $G$. Similarly, the \emph{stable torsion length} in $G$ is the stable length with respect to the set $T$ of torsion elements in $G$.
The stable $S$-length of an element $g \in G$ is strictly positive, as soon as there exists a homomorphism $\phi: G \to \Z$, bounded on $S$, with $\phi(g) \neq 0$. Many groups do not support any non-trivial integer-valued homomorphisms, e.g. $SL(2,\Z)$. However, some interesting groups support real-valued quasimorphisms, which provide lower bounds for stable lengths, as well.
In \cite{Ko}, Kotschick proved several strong results on the mapping class groups of closed surfaces, including the strict positivity of the stable torsion length of Dehn twists. These results are based upon the existence of quasi-morphisms. In the remaining three sections, we study stable commutator and torsion lenghts via the Rasmussen invariant of braids.

\begin{theorem} The Rasmussen invariant $s$ is a quasi-morphism of defect $n+1$ on the braid group $B_n$:
$$|s(\alpha \beta)-s(\alpha)-s(\beta)| \leq n+1,$$
for all $\alpha$, $\beta \in B_n$.
\end{theorem}

We apply a combination of the Rasmussen invariant and the signature of braids to estimate the stable commutator and torsion lengths in braid groups modulo their centers $B_n/ \langle \Delta_n \rangle$. The center of the braid group $B_n$ is generated by the element 
$\Delta_n=(\sigma_1 \sigma_2 \ldots \sigma_{n-1})^n$ (a `full twist'). For an element $w \in B_n$, we denote its image in $B_n/ \langle \Delta_n \rangle$ under the natural projection by $\overline{w}$.

\begin{proposition} Let $n \geq 2$.
\begin{enumerate}
\item For any element of the braid group $w \in B_n$, $w^{n(n-1)}$ is a product of commutators in $B_n/ \langle \Delta_n \rangle$. \\
\item $B_n/ \langle \Delta_n \rangle$ is generated by the two torsion elements
$\overline{\sigma_1 \sigma_2 \ldots \sigma_{n-1}}$ 
and $\overline{\sigma_1 \sigma_2 \ldots \sigma_{n-1} \sigma_1}$. \\
\end{enumerate}
\end{proposition}

A proof of this proposition is presented in the appendix. 

\begin{remark} The braid group itself has no torsion elements, since it admits a left-invariant order.
\end{remark}

Proposition~1 allows us to speak of the stable commutator length 
$$c: B_n/ \langle \Delta_n \rangle \to \R,$$
and the stable torsion length 
$$t: B_n/ \langle \Delta_n \rangle \to \R.$$ 
We shall estimate the stable commutator and torsion lengths for alternating braids. Here a braid is \emph{alternating}, if it can be written as a product of `alternating' generators $\sigma_1$, $\sigma_2^{-1}$,  
$\sigma_3$, $\sigma_4^{-1}$, $\ldots$, or $\sigma_1^{-1}$, $\sigma_2$, $\sigma_3^{-1}$, $\sigma_4$, $\ldots$.

\begin{theorem} Let $\alpha \in B_n$ be an alternating braid with $lk(\alpha) \neq 0$, $n \geq 3$. Then the stable commutator and torsion lengths of $\overline{\alpha} \in B_n/ \langle \Delta_n \rangle$ are strictly positive.
\end{theorem}

The assumption $lk(\alpha) \neq 0$ is absolutely essential:
in the braid groups on more than 4 strings there exist alternating braids which are conjugate to their inverses, for example, the braid $\alpha=\sigma_1 \sigma_4^{-1}$. Indeed,
$$\alpha^{-1}=\beta \alpha \beta^{-1},$$
for $\beta=\sigma_2 \sigma_3 \sigma_4 \sigma_1 \sigma_2 \sigma_3 \sigma_1^{-1} \sigma_2^{-1}$, see figure~3.
Therefore, 
$$\alpha^{2l}=\alpha^l \alpha^l=\alpha^l \beta^{-1} \alpha^{-l} \beta$$ 
is a commutator, for all $l \in \N$. This implies $c(\alpha)=0$.

\begin{figure}[ht]
\scalebox{1}{\raisebox{-0pt}{$\vcenter{\hbox{\epsffile{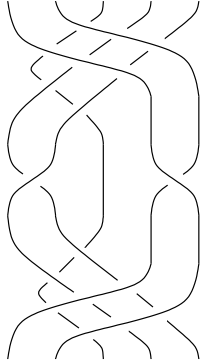}}}$}} 
\caption{$\beta \alpha \beta^{-1}=\alpha^{-1}$}
\end{figure}

The case $n=3$ has been well-studied, since $B_3/ \langle \Delta_3 \rangle$ is isomorphic to $PSL(2, \Z)$, via the correspondence
$$
\overline{\sigma}_1 \mapsto
\begin {pmatrix}
 1 & 1 \\
 0 & 1 \\
\end{pmatrix},
$$
$$
\overline{\sigma}_2 \mapsto
\begin {pmatrix}
 1 & 0 \\
 -1 & 1 \\
\end{pmatrix}.
$$
In \cite{PR}, Polterovich and Rudnick proved much stronger statements than our theorem~2 for $SL(2, \Z)$. For instance, every non-torsion element of $SL(2, \Z)$, which is not conjugate to its inverse, can be detected by a homogeneous real-valued quasi-morphism. This implies the strict positivity of the stable commutator and torsion lengths of these elements (see~\cite{Ko} for a detailed background on quasi-morphisms and stable lengths).

By a recent result of Endo and Kotschick (\cite{EK}), there exist non-torsion elements in the mapping class groups of closed hyperbolic surfaces, which are not conjugate to their inverses and still have vanishing stable commutator length. As they explain, whenever three commuting elements $a$, $b$, $c$ of an arbitrary group are pairwise conjugate, then the element $abc^{-2}$ has vanishing stable commutator length. We can easily carry over their example to braid groups on more than 5 strings.

\begin{example} The braid $\sigma_1 \sigma_3 \sigma_5^{-2} \in B_6$ has infinite order and vanishing stable commutator length. Moreover, two different powers of $\sigma_1 \sigma_3 \sigma_5^{-2}$ are never conjugate in $B_6$, since the natural closure of $(\sigma_1 \sigma_3 \sigma_5^{-2})^m$ is a 3- or 6-component link consisting of two $(2,m)$-torus links and one $(2,-2m)$-torus link.
\end{example}

By the slice-Bennequin inequality (\cite{Ru}) and Beliakova-Wehrli's results on the Rasmussen invariant (\cite{BW}), we know that
$$s(L(\beta))=1+lk(\beta)-n,$$
for all positive braids $\beta \in B_n$. Here $L(\beta)$ denotes the natural closure of $\beta$ in $\R^3$. Hence, all positive braids are mapped to zero by the quasi-morphism $s-lk+n-1$ (since $\Delta_n$ is positive). In particular, 
$s-lk+n-1$ descends to a quasi-morphism on $B_n/ \langle \Delta_n \rangle$. This raises the following questions.

\begin{question} Is $s-lk+n-1$ bounded or not on $B_n/ \langle \Delta_n \rangle$?
\end{question}

\begin{question} What is the interpretation of $s-lk+n-1$ on $PSL(2, \Z) \simeq B_3/ \langle \Delta_3 \rangle$?
\end{question}

\section{Rasmussen invariant and signature as quasi-morphisms}

Recently, Beliakova and Wehrli extended the definition of the Rasmussen invariant $s$ to links (\cite{BW}). As in the case of knots, $s$ provides a lower bound for the 4-genus of links, as follows: let $S$ be a smooth oriented cobordism between two oriented links $L_0$,$L_1 \subset \R^3$, i.e. a smooth oriented properly embedded surface $S \subset \R^3 \times [0,1]$, such that the oriented boundary components $S \cap \R^3 \times \{i\}$ correspond to the oriented links $L_i$. Then
\begin{equation}
|s(L_1)-s(L_0)| \leq -\chi(S),
\label{genusbound}
\end{equation}
where $\chi(S)$ is the Euler characteristic of $S$.
We define the Rasmussen invariant for braids by
$$s(\alpha)=s(L(\alpha)).$$
As before, $L(\alpha)$ denotes the natural closure of $\alpha$ in $\R^3$.
Inequality~(\ref{genusbound}) allows us to prove that $s$ is a quasimorphism on $B_n$. Indeed, for any two braids $\alpha$, $\beta \in B_n$, there is a cobordism $S(\alpha, \beta)$ between $L(\alpha \beta)$ and the disjoint union of $L(\alpha)$ and 
$L(\beta)$. $S(\alpha, \beta)$ consists of a disjoint union of annuli, glued together by $n$ bands. Figure~4 shows the cobordism $S(\alpha, \beta)$ for $\alpha=\sigma_1 \sigma_2^{-1} \sigma_1$ and $\beta=\sigma_2 \sigma_1$. It consists of three annuli, glued together by three bands located in the middle of the figure.
(In fact, $S(\alpha, \beta)$ is a cobordism between $L(\alpha \beta)$ with the opposite orientation and the disjoint union of $L(\alpha)$ and $L(\beta)$. However, this does not matter, since a global change of orientation for a link does not affect its Rasmussen invariant.)
The Euler characteristic of $S(\alpha, \beta)$ equals $-n$. Therefore,
$$|s(L(\alpha \beta))-s(L(\alpha) \amalg L(\beta))| \leq n.$$
Here $L(\alpha) \amalg L(\beta)$ denotes the disjoint union of the two links $L(\alpha)$ and $L(\beta)$. By Beliakova and Wehrli's definition (\cite{BW}),
$$s(L(\alpha) \amalg L(\beta))=s(L(\alpha))+s(L(\beta))-1.$$
This proves theorem~1.

\begin{figure}[ht]
\scalebox{1}{\raisebox{-0pt}{$\vcenter{\hbox{\epsffile{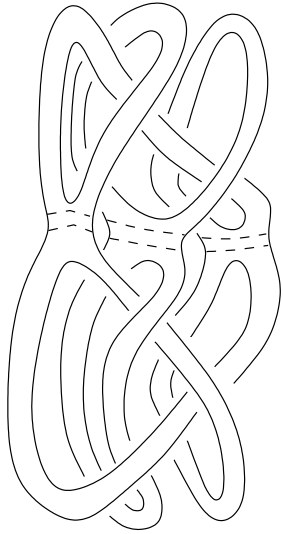}}}$}} 
\caption{cobordism $S(\alpha, \beta)$}
\end{figure}

The signature is an integer-valued link invariant. It is defined as the signature of any symmetrized Seifert matrix of a given link (see~\cite{Ka}). We may define a signature for braids, as well, via their natural closure in $\R^3$. In \cite{GG2}, Gambaudo and Ghys proved that the signature on $B_n$ is a quasi-morphism of defect $n$, i.e. 
$$|\sigma(\alpha \beta)-\sigma(\alpha)-\sigma(\beta)| \leq n,$$
for all $\alpha$, $\beta \in B_n$.
The signature of a trivial link equals zero. A computation of the signature of torus links was done by Gordon, Litherland and Murasugi in \cite{GLM}. Their result implies the following asymptotical behaviour of the signature on powers of the central element $\Delta_n$, $n \geq 2$:
\begin{equation}
\lim_{l \to \infty} \frac{1}{l} \; \sigma(\Delta_n^l)=\frac{1}{2} (n^2-1), \text{ if n is odd,}
\label{oddsign}
\end{equation}
\begin{equation}
\lim_{l \to \infty} \frac{1}{l} \; \sigma(\Delta_n^l)=\frac{1}{2} n^2, \text{ if n is even.}
\label{evensign}
\end{equation}
We shall use these two equalities in the proof of theorem~2.

\section{Quasi-morphisms and stable lengths}

A quasi-morphism $\phi: G \to \R$ is \emph{homogeneous}, if $\phi(g^n)=n \phi(g)$, for all $n \in \N$, $g \in G$. 
As explained in \cite{Ko}, the stable torsion length $t(g)$ of a group element $g \in G$ is strictly positive, as soon as there exists a homogeneous quasi-morphism $\phi: G \to \R$ with $\phi(g) \neq 0$. More precisely, 
\begin{equation}
c(g)>\frac{|\phi(g)|}{2D(\phi)},
\label{commutator}
\end{equation}
\begin{equation}
t(g)>\frac{|\phi(g)|}{D(\phi)},
\label{torsion}
\end{equation}
where $D(\phi)$ is the defect of $\phi$.
Every quasi-morphism $\phi: G \to \R$ can be homogenized by the following procedure:
$$\widetilde{\phi} (g)=\lim_{n \to \infty} \frac{1}{n} \phi(g^n).$$
We denote the homogenized Rasmussen invariant and signature by $\widetilde{s}$ and $\widetilde{\sigma}$, respectively. By the slice-Bennequin inequality,
$$s(\Delta_n^l)=1+lk(\Delta_n^l)-n=1+ln(n-1)-n,$$
whence
\begin{equation}
\widetilde{s}(\Delta_n)=n(n-1).
\label{stilde}
\end{equation}
In case $n$ is odd, $\widetilde{\sigma}(\Delta_n)=\frac{1}{2}(n^2-1)$, by~(\ref{oddsign}).\\
In case $n$ is even, $\widetilde{\sigma}(\Delta_n)=\frac{1}{2}n^2$, by~(\ref{evensign}).\\
In the former case,
\begin{equation}
\phi_{odd}:=2n(n-1) \widetilde{\sigma}-(n^2-1) \widetilde{s}
\label{phiodd}
\end{equation}
descends to a homogeneous quasi-morphism on $B_n/ \langle \Delta_n \rangle$.\\
In the latter case,
\begin{equation}
\phi_{even}:=2n(n-1) \widetilde{\sigma}-n^2 \widetilde{s}
\label{phieven}
\end{equation}
descends to a homogeneous quasi-morphism on $B_n/ \langle \Delta_n \rangle$.\\

\begin{proof}[Proof of theorem~2] Let $\alpha \in B_n$, $n \geq 3$, be an alternating braid with $lk(\alpha) \neq 0$. 
We arrange $lk(\alpha)>0$, by taking the mirror image of $\alpha$, if necessary.
We claim that $\widetilde{s}(\alpha)=\widetilde{\sigma}(\alpha)$. Indeed, $\alpha^l$ is alternating, for all $l \in \N$. If 
$L(\alpha^l)$ is a knot, then $s(\alpha^l)=\sigma(\alpha^l)$, by property~(2) of the Rasmussen invariant. In case 
$L(\alpha^l)$ has several components, we may multiply $\alpha^l$ by a suitable alternating braid $\beta \in B_n$ of length 
$\leq n-1$, such that $L(\alpha^l \beta)$ is an alternating knot. Hence, the asymptotical behaviour of the Rasmussen invariant and the signature on powers of $\alpha$ are the same:
$$\widetilde{s}(\alpha)=\widetilde{\sigma}(\alpha).$$
The expressions for $\phi_{odd}(\alpha)$ and $\phi_{even}(\alpha)$ simplify to
$$\phi_{odd}(\alpha)=(n^2-2n+1) \widetilde{s}(\alpha),$$
in case $n$ is odd;
$$\phi_{even}(\alpha)=(n^2-2n) \widetilde{s}(\alpha),$$
in case $n$ is even. The additional assumption $lk(\alpha)>0$ ensures $\widetilde{s}(\alpha)>0$. Indeed,
$$s(\alpha^l) \geq 1+lk(\alpha^l)-n=1+lk(\alpha)l-n,$$
whence $\widetilde{s}(\alpha) \geq lk(\alpha)>0$. At last, if 
$n \geq 3$, then $\phi_{odd}(\alpha)>0$ or $\phi_{even}(\alpha)>0$. This proves theorem~2.
\end{proof}

\begin{remark} 
The defects of the homogeneous quasi-morphisms $\phi_{odd}$ and $\phi_{even}$ on $B_n/ \langle \Delta_n \rangle$ are proportional to $n^3$, since the defects of $\widetilde{\sigma}$ and $\widetilde{s}$ are proportional to $n$. Using the estimates (\ref{commutator}) and (\ref{torsion}), we conclude that the stable commutator and torsion lengths of an alternating braid $\overline{\alpha} \in B_n/ \langle \Delta_n \rangle$ are bounded below by a constant times $\frac{1}{n}$. In particular, this bound tends to zero, if $n$ tends to infinity. We point out that Kotschick's lower bound for the stable torsion length of products of right-handed Dehn twists in the mapping class groups shows an analogous behaviour, with respect to the genus.
\end{remark}

\begin{appendix}
\section{The braid group modulo its center}
It is well-known that the braid group $B_n$ is generated by the two elements $\sigma_1 \sigma_2 \ldots \sigma_{n-1}$ and 
$\sigma_1$. Therefore, $B_n/ \langle \Delta_n \rangle$ is generated by the two elements $a=\overline{\sigma_1 \sigma_2 \ldots \sigma_{n-1}}$ and 
$b=\overline{\sigma_1 \sigma_2 \ldots \sigma_{n-1} \sigma_1}$. We observe that the order of $a \in B_n/ \langle \Delta_n \rangle$ is $n$:
$$a^n=(\sigma_1 \sigma_2 \ldots \sigma_{n-1})^n=\Delta_n.$$ 
Further, the order of $b \in B_n/ \langle \Delta_n \rangle$ is $n-1$:
$$b^{n-1}=(\sigma_1 \sigma_2 \ldots \sigma_{n-1} \sigma_1)^{n-1}=\Delta_n.$$ 
The latter equality is easy to see, as shown in figure~5 for $n=4$ (the two braids $b^{n-1}$ and $\Delta_n$ differ by an isotopy of the second strand). This proves the second statement of proposition~1.

\begin{figure}[ht]
\scalebox{1}{\raisebox{-0pt}{$\vcenter{\hbox{\epsffile{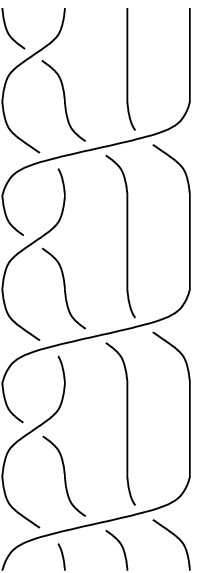}}}$}} \qquad $=$ \qquad
\scalebox{1}{\raisebox{-0pt}{$\vcenter{\hbox{\epsffile{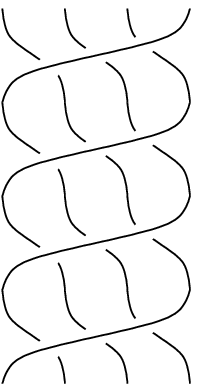}}}$}}
\caption{$(\sigma_1 \sigma_2 \sigma_3 \sigma_1)^3=\Delta_4$}
\end{figure}

As to the first statement of proposition~1, we observe that all the elements $\sigma_i^m \sigma_j^{-m}$ are commutators, for all $i,j \in \{1,2,\ldots,n-1\}$ (this easily follows from the fact that $\sigma_i$ and $\sigma_j$ are conjugate in $B_n$).
The elements $\sigma_i^m \sigma_j^{-m}$ generate the kernel of the linking number $lk: B_n \to \Z$, i.e. $\ker(lk)=[B_n,B_n]$. Let $w \in B_n$ be any element of the braid group $B_n$. The linking number of $w^{n(n-1)}$ is a multiple of $n(n-1)=lk(\Delta_n)$, since 
$lk(w^{n(n-1)})=n(n-1) lk(w)$. Therefore, we may write 
$$w^{n(n-1)}=\Delta_n^{lk(w)} \beta,$$ 
where $\beta$ is a braid with $lk(\beta)=0$, i.e. $\beta \in [B_n,B_n]$. This completes the proof of proposition~1.
\end{appendix}

\section*{Acknowledgements}
I would like to thank \'Etienne Ghys and Dieter Kotschick for the fruitful conversations we have had on asymptotic invariants, mapping class groups and braid groups.

\bigskip
\noindent
Department of Mathematics,
ETH Z\"urich, 
Switzerland

\bigskip
\noindent
\emph{sebastian.baader@math.ethz.ch}

\end{document}